\documentclass[letterpaper, 10 pt, conference]{ieeeconf}  

\IEEEoverridecommandlockouts                              

\title{\LARGE \bf
GPS-Denied Relative Motion Estimation For Fixed-Wing UAV
Using the Variational Pose Estimator
}

\author{Maziar Izadi$^{1}$, Amit K. Sanyal$^{2,\dag}$, Randy Beard$^{3}$ and He Bai$^{4}$
\thanks{$^{1}$M. Izadi is with the Department of Mechanical and Aerospace Engineering, New Mexico State University, Las Cruces, NM 88003 USA.
        {\tt\small mi@nmsu.edu}}%
\thanks{$^{2}$A. K. Sanyal is with the Department of Mechanical and Aerospace Engineering, Syracuse University, Syracuse, NY 13244 USA.
        {\tt\small aksanyal@syr.edu}}%
\thanks{$^{3}$R. Beard is with the Electrical and Computer Engineering Department, Brigham Young University, Provo, UT 84602 USA.
        {\tt\small beard@ee.byu.edu}}%
\thanks{$^{4}$H. Bai is with UtopiaCompression Corporation, Los Angeles, CA 90064 USA. 
        {\tt\small he@utopiacompression.com}} 
\thanks{$^{\dag}$Address all correspondence to this author.}
}

\usepackage{graphicx}          
\usepackage{amsmath,amsfonts,amssymb}
\usepackage{mathtools, textcomp}
\usepackage{mathrsfs}
\usepackage{graphicx}
\usepackage[dvips]{epsfig}
\usepackage{epstopdf}
\usepackage{color}
\usepackage{dsfont}

\newcommand{\SO}{\ensuremath{\mathrm{SO(3)}}}

\newcommand{\Ta}{\ensuremath{\mathrm{T}}}
\newcommand{\T}{^{\mbox{\small T}}}
\newcommand{\so}{\ensuremath{\mathfrak{so}(3)}}
\newcommand{\SE}{\ensuremath{\mathrm{SE(3)}}}

\newcommand{\bR}{\ensuremath{\mathbb{R}}}
\newcommand{\bS}{\ensuremath{\mathbb{S}}}

\newcommand{\mrm}{\mathrm}

\newcommand{\diag}{\mbox{diag}}

\newcommand{\bbm}{\begin{bmatrix}}
\newcommand{\ebm}{\end{bmatrix}}
\newcommand{\matl}{\left[ \begin{array}}
\newcommand{\matr}{\end{array} \right]}
\newcommand{\be}{\begin{equation}}
\newcommand{\ee}{\end{equation}}
\newcommand{\bea}{\begin{eqnarray}}
\newcommand{\eea}{\end{eqnarray}}
\newcommand{\beas}{\begin{eqnarray*}}
\newcommand{\eeas}{\end{eqnarray*}}
\newcommand{\nn}{\nonumber}
\newcommand{\mC}{\mathcal{C}}

\newcommand{\cL}{\mathcal{L}}

\newcommand{\cT}{\mathcal{T}}
\newcommand{\cU}{\mathcal{U}}

\newcommand{\di}{\mathrm{d}}

\newcommand{\tr}{\mathrm{trace}}

\newcommand{\lan}{\langle}
\newcommand{\ran}{\rangle}

\newcommand{\cS}{\mathcal{S}}

\newcommand{\ad}[1]{{\mathrm{ad}_{#1}}}          			
\newcommand{\adast}[1]{{\mathrm{ad}_{#1}^\ast}}  			
\newcommand{\Ad}[1]{{\mathrm{Ad}_{#1}}}  			
\DeclareMathOperator{\expm}{{expm}}  			

\newtheorem{theorem}{Theorem}[section]

\def\thesection{\arabic{section}}

\newcommand{\bi}{\begin{itemize}}
\newcommand{\ei}{\end{itemize}}

\DeclareMathAlphabet{\mathpzc}{OT1}{pzc}{m}{it}

\newcommand{\bJ}{\mathbb{J}}

\newcommand{\mpz}{\mathpzc}
\newcommand{\msg}{\mathsf{g}}
\newcommand{\msh}{\mathsf{h}}

\newcommand{\sS}{\mathsf{S}}
\newcommand{\sO}{\mathsf{O}}
\newcommand{\bD}{\mathbb{D}}
\begin{document}

\maketitle
\thispagestyle{empty}
\pagestyle{empty}

\begin{abstract}

Relative pose estimation between fixed-wing unmanned aerial vehicles (UAVs) is treated 
using a stable and robust estimation scheme. The motivating application of this scheme is that 
of ``handoff" of an object being tracked from one fixed-wing UAV to another in a team of UAVs, 
using onboard sensors in a GPS-denied environment. This estimation scheme uses optical 
measurements from cameras onboard a vehicle, to estimate both the relative pose and relative 
velocities of another vehicle or target object. It is obtained by applying the Lagrange-d'Alembert 
principle to a Lagrangian constructed from measurement residuals using only the optical 
measurements. This nonlinear pose estimation scheme is discretized for computer implementation 
using the discrete Lagrange-d'Alembert principle, with a discrete-time linear filter for obtaining 
relative velocity estimates from optical measurements. Computer simulations depict the stability 
and robustness of this estimator to noisy measurements and uncertainties in initial relative pose 
and velocities.

\end{abstract}

\section{INTRODUCTION}\label{Sec1}
Onboard estimation of relative motion between unmanned vehicles and spacecraft is an 
important enabling technology for autonomous operations of teams and formations of such 
vehicles. A stable relative motion estimation scheme that is robust to measurement noise and 
requires no knowledge of the dynamics model of the vehicle being observed, is presented here. 
This estimation scheme can enhance the autonomy and reliability of teams of unmanned 
vehicles operating in uncertain GPS-denied environments. Salient features of this estimation 
scheme are: (1) use of only onboard optical sensors for estimation of relative pose and 
velocities; (2) robustness to uncertainties and lack of knowledge of dynamics model of 
observed vehicle; (3) low computational complexity such that it can be implemented 
with onboard processors; and (4) proven stability with large domain of attraction for relative 
motion state estimation errors. Stable and robust relative motion estimation of unconstrained 
motion of teams of unmanned vehicles in the absence of complete knowledge of their 
dynamics, is required for their safe, reliable, and autonomous operations in poorly known 
environments. In practice, the dynamics of an observed vehicle may not be perfectly known, 
especially in outdoor environments where the vehicle may be under the action of unknown 
forces and moments. The scheme proposed here has a single, stable algorithm for the 
naturally coupled relative translational and rotational motion between unmanned vehicles, 
using measurements from onboard optical sensors. This avoids the need for 
measurements from external sources, like GPS, which may not be available in indoor, 
underwater or cluttered environments \cite{amelin2014algorithm,leishman2014relative,miller2014tracking}.

Relative pose (position and attitude) estimation of one vehicle from another vehicle is treated 
here. Determining the relative attitude requires that at least three feature points on the observed 
vehicle are available. Attitude estimation and control schemes that use generalized coordinates 
or quaternions for attitude representation are usually {\em unstable in the sense of Lyapunov}, 
as has been shown in recent research~\cite{Bayadi2014almost,chaturvedi2011rigid,sanyal2009inertia}. 
One adverse consequence of these unstable estimation and control schemes is that they end 
up taking longer to converge compared to stable schemes with the same initial conditions and 
same initial transient behavior. Attitude observers and filtering schemes on $\SO$ and $\SE$ 
have been reported in, 
e.g., \cite{bonmaro09,Khosravian2015Recursive,Khosravian2015observers,mahapf08,mabeda04,markSO3,rehbinder2003pose,sany,Vas1,silvest08}. 
These estimators do not suffer from kinematic singularities like estimators using coordinate 
descriptions of attitude, and they do not suffer from the unstable unwinding phenomenon 
encountered by continuous estimators using unit quaternions. Recently, the 
maximum-likelihood (minimum energy) filtering method of Mortensen~\cite{Mortensen} was 
applied to attitude and pose estimation on $\SO$ and $\SE$, resulting in nonlinear estimation 
schemes that seek to minimize the stored ``energy" in measurement errors 
\cite{aguhes06,ZamPhD,Zamani2013minimum}. This led to ``near optimal"  filtering schemes 
that are based on approximate solutions of the Hamilton-Jacobi-Bellman (HJB) equation and 
do not have provable stability. The estimation scheme obtained here is shown to be almost 
globally asymptotically stable. Moreover, unlike filters based on Kalman filtering, the 
estimator proposed here does not make any assumptions on the statistics of initial state 
estimate or sensor noise. 

For the relative pose estimation problem analyzed in this paper, it is assumed that one vehicle 
can optically measure a known pattern fixed to the body of another vehicle whose relative 
motion states are to be estimated. From such optical (camera) measurements, the relative 
velocities (translational and angular) are also estimated. The variational attitude estimator 
recently appeared \cite{ACC2015,Automatica}, where it was shown to be almost globally 
asymptotically stable. The advantages of this scheme over Kalman-based schemes are 
reported in \cite{ICRA2015}. A companion paper extends the variational attitude estimator to 
estimation of coupled rotational (attitude) and translational motion. Maneuvering vehicles, like 
UAVs tracking ground targets, have naturally coupled rotational and translational
motion. In such applications, designing separate state estimators for the translational and 
rotational motions may not be effective and could lead to poor navigation. For relative pose 
estimation between such vehicles operating in teams, the approach proposed here for robust 
and stable estimation will be more effective than Kalman filtering-based schemes. 
The estimation scheme proposed here can be implemented without any velocity measurements,  
which is useful when Doppler lidar sensors are not available onboard or rate gyros are corrupted 
by high noise content and bias \cite{Good2013ECC,Good2014ASME,Good2014ACC}. 


\section{RELATIVE NAVIGATION USING OPTICAL SENSORS} \label{Sec2}
\subsection{Motivation}
When multiple UAV perform surveillance and target tracking missions in a GPS-denied 
environment, they need to ensure that they are tracking the same target of interest. When 
necessary, the tracking responsibility may need to be handed off from one UAV to another. 
When GPS signals are available, such a handoff procedure can be achieved by a \textit{tracking 
UAV} geo-locating the target and sending the global coordinates of the target to a \textit{handoff 
UAV}. In GPS-denied environments, such a handoff procedure faces several challenges. The 
most significant challenge is the following: because no GPS signals are available, the 
\textit{handoff UAV} may not know the position of the \textit{tracking UAV}. Therefore, it needs 
to use on-board sensors to detect and navigate towards the \textit{tracking UAV}. Moreover, 
global information about the target is not available. Because the \textit{tracking UAV} does not 
have GPS, it can only geo-locate the target in its own navigation frame. Since the 
\textit{handoff UAV} has a different coordinate system than the \textit{tracking UAV}, the target 
information from the \textit{tracking UAV} cannot be directly used by the \textit{handoff UAV} to 
track the target. The \textit{handoff UAV} has to perform a coordinate transformation that 
converts the target information to its own navigation frame. This task is carried out by the 
relative pose estimation technique presented here. 


\subsection{Relative Pose Measurement Model} 
Let $O$ denote the observed vehicle and $S$ denote the vehicle that is observing $O$. Let 
$\mathsf{S}$ denote a coordinate frame fixed to $S$ and $\sO$ be a coordinate frame fixed to 
$O$. Let $R\in\SO$ be the rotation matrix from frame $\sS$ to frame $\sO$ and $b$ denote 
the position of origin of $\sS$ expressed in frame $\sO$. The pose (transformation) of frame 
$\sS$ to frame $\sO$ is 
\begin{align} 
&\msg=\bbm R & b\\ 0 & 1\ebm\in\SE.
\label{gDef} 
\end{align}
The positions of a fixed set of feature points or patterns on vehicle $O$ are observed by optical 
sensors fixed to vehicle $S$. Velocities of these points are not directly measured, but may 
be calculated using a simple linear filter as in \cite{ACC2015}. Assume that there are $\mpz j>2$ 
feature points, which are always in the sensor field-of-view (FOV) of the sensor fixed to vehicle 
$S$, and the positions of these points are known in frame $\sO$ as $p_j$, $j\in\{1,2,\ldots,\mpz j\}$. 
These points generate ${\mpz j\choose 2}$ unique pairwise relative position vectors, which are the 
vectors connecting any two of these points. 

Denote the position of the optical sensor on vehicle $S$ and the vector from that sensor to 
an observed point on vehicle $O$ as $s\in\bR^3$ and $q_j\in\bS^2$, $j=1,2,\ldots,\mpz j$, 
respectively, both vectors expressed in frame $\sS$. Thus, in the absence of measurement 
noise
\begin{align}
p_j=R(q_j+s)+b=Ra_j+b,\; j\in\{1,2,\ldots,\mpz j\},
\label{FrameTrans}
\end{align}
where $a_j=q_j+s$, are positions of these points expressed in $\sS$. In practice, the $a_j$ are 
obtained from proximity optical measurements that will have additive noise; denote by 
$a_j^m$ the measured vectors. The mean values of the vectors $p_j$ and $a_j^m$ are denoted 
as $\bar{p}$ and $\bar{a}^m$, and satisfy
\begin{align}
\bar{a}^m= R\T (\bar p- b)+ \bar\varsigma, \label{barp}
\end{align}
where $\bar{p}=\frac{1}{\mpz j}\sum\limits^\mpz j_{j=1}p_j$,  $\bar{a}^m=\frac{1}{\mpz j}\sum\limits^\mpz j_{j=1}
a_j^m$  and $\bar\varsigma$ is the additive measurement noise 
obtained by averaging the measurement noise vectors for each of the $a_j$. Consider the ${\mpz j\choose 2}$ relative position vectors from optical measurements, 
denoted as $d_j=p_\lambda-p_\ell$ in frame $\sO$ and the corresponding vectors in frame $\sS$ as 
$l_j=a_\lambda-a_\ell$, for $\lambda,\ell\in\{1,2,\ldots,\mpz j\}$, $\lambda\ne \ell$. Therefore,
\begin{align}
d_j=Rl_j\Rightarrow D=RL,
\end{align}
where $D=[d_1\;\, \cdots\;\, d_\beta]$, $L=[l_1\;\, \cdots\;\, l_\beta]\in\bR^{3\times\beta}$ with $\beta={\mpz j\choose 2}$. Note that the matrix of known relative vectors $D$ is assumed to be known 
and bounded. Denote the measured value of matrix $L$ in the presence of measurement noise 
as $L^m$. Then,
\begin{align}
L^m=R\T D+\mathscr{L},
\label{VecMeasMod}
\end{align}
where $\mathscr{L}\in\bR^{3\times\beta}$ is the matrix of measurement errors in these vectors observed 
in frame $\sS$. 

\subsection{Relative Velocities Measurement Model} 
Denote the relative angular and translational velocity of vehicle $O$ expressed in frame $\sS$ by $\Omega$ and $\nu$, respectively. Thus, one can write the kinematics 
of the rigid body as
\begin{align}
\dot{\Omega}=R\Omega^\times,\dot{b}=R\nu\Rightarrow\dot{\msg}= \msg\xi^\vee,
\label{Kinematics}
\end{align}
where $\xi= \bbm \Omega\\ \nu\ebm\in\bR^6$ and $\xi^\vee=\bbm \Omega^\times &\; \nu\\ 0 \;\;& 0\ebm$ and $(\cdot)^\times: \bR^3\to\so\subset\bR^{3\times 3}$ is the 
skew-symmetric cross-product operator that gives the vector space isomorphism between 
$\bR^3$ and $\so$. In order to do 
so, one can differentiate \eqref{FrameTrans} as follows
\begin{align}
&\dot{p}_j=R\Omega^\times a_j+R\dot{a}_j+\dot{b}=R\big(\Omega^\times a_j+\dot{a}_j+\nu\big)=0\nn\\
\Rightarrow&\dot{a}_j-a_j^\times\Omega+\nu=0\nn\\
\Rightarrow&v_j=\dot{a}_j=[a_j^\times\; -I]\xi=G(a_j)\xi,
\label{pjdot}
\end{align}
where $G(a_j)=[a_j^\times\; -I]$ has full row rank. From vision-based or Doppler lidar sensors, one can also 
measure the velocities of the observed points in frame $\sS$, denoted $v_i^m$. Here, velocity 
measurements as would be obtained from vision-based sensors is considered. 
The measurement model for the velocity is of the form
\be v_j^m=G(a_j)\xi+\vartheta_j,\ee
where $\vartheta_j\in\bR^3$ is the additive error in velocity measurement $v_j^m$. Instantaneous 
angular and translational velocity determination from such measurements is treated 
in~\cite{ast_acc14}. Note that $v_j=\dot{a}_j$, for $j\in\{1,2,\ldots,\mpz j\}$. As this kinematics indicates, the relative 
velocities of at least three beacons are needed to determine the vehicle's translational and 
angular velocities uniquely at each instant. The rigid body velocities are obtained 
using the pseudo-inverse of $\mathds{G}(A^f)$:
\begin{align}
\mathds{G}(A^f)\xi^f&=\mathds{V}(V^f)\Rightarrow\xi^f=\mathds{G}^\ddag(A^f)\mathds{V}(V^f),\label{ximMore2}\\
\mbox{where }\;\mathds{G}(A^f)&=\bbm G(a^f_1)\\\vdots\\G(a^f_\mpz j)\ebm\;\mbox{and }\;\mathds{V}(V^f)=\bbm v^f_1\\\vdots\\v^f_\mpz j \ebm.\label{GVDef}
\end{align}
When at least three beacons are measured, $\mathds{G}(A^f)$ is a full column rank matrix, 
and $\mathds{G}^\ddag(A^f)= \Big( \mathds{G}\T(A^f)\mathds{G} (A^f)\Big)^{-1} 
\mathds{G}\T (A^f)$ gives its pseudo-inverse. For the case that only one or two beacons are 
observed, $\mathds{G}(A^f)$ is a full row rank matrix, whose pseudo-inverse is given by 
$\mathds{G}^\ddag (A^f)= \mathds{G}\T (A^f)\Big( \mathds{G}(A^f)\mathds{G}\T 
(A^f)\Big)^{-1}$.

\section{DYNAMIC ESTIMATION OF MOTION FROM PROXIMITY MEASUREMENTS}\label{Sec3}
In order to obtain state estimation schemes from measurements as outlined in Section 
\ref{Sec2} in continuous time, the Lagrange-d'Alembert principle is applied to an action functional of 
a Lagrangian of the state estimate errors, with a dissipation term linear 
in the velocities estimate error. This section presents the estimation scheme obtained 
using this approach. Denote the estimated pose and its kinematics as
\begin{align}
\hat\msg=\bbm \hat R & \;\;\; \hat{b}\\ 0 & \;\;\; 1\ebm\in\SE, \;\;\dot{\hat\msg}=\hat\msg\hat\xi^\vee,
\end{align}
where $\hat\xi$ is rigid body velocities estimate, with $\hat\msg_0$ as the initial pose estimate and the pose estimation 
error as
\begin{align}
\msh=\msg\hat\msg^{-1}=\bbm Q &\;\;\;\; b-Q\hat{b}\\ 0 & \;\;1\ebm=\bbm Q &\;\;\; x\\ 0 &\;\;\;1\ebm\in\SE,
\end{align}
where $Q=R\hat{R}\T$ is the attitude estimation error and $x=b-Q\hat{b}$. Then one obtains, in the case of perfect measurements,
\begin{align}\begin{split}
&\dot\msh = \msh\varphi^\vee,\, \mbox{ where }\, \varphi(\hat\msg,\xi^m,\hat\xi)=\bbm\omega\\\upsilon\ebm= \Ad{\hat\msg}\big(\xi^m- \hat\xi),
\end{split}\label{hdot}
\end{align}
where $\Ad{\mpz{g}}=\bbm \mpz{R}~&~0\\ \mpz{b}^\times\mpz{R}~&~\mpz{R}\ebm$ for $\mpz{g}=\bbm \mpz{R} &\;\; \mpz{b}\\ 0 & \;\;1\ebm$. The attitude and position estimation error dynamics are also in the form
\begin{align}
\dot{Q}=Q\omega^\times,\;\;\dot{x}=Q\upsilon.
\label{Qdot}
\end{align}

\subsection{Lagrangian from Measurement Residuals}
Consider the sum of rotational and translational measurement residuals between the 
measurements and estimated pose as a potential energy-like function. Defining the trace inner 
product on $\mathbb{R}^{n_1\times n_2}$ as
\begin{align}
\lan A_1,A_2\ran :=\tr(A_1 \T A_2),\label{tr_def}
\end{align}
the rotational potential function (Wahba's cost function~\cite{jo:wahba}) is expressed as
\begin{align}
\cU^0_r (\hat{\msg},L^m,D) &= \frac12\lan D -\hat R L^m , (D -\hat R L^m)W\ran,
\label{U0r}
\end{align}
where $W=\diag(w_j)\in\bR^{n\times n}$ is a positive diagonal matrix of weight factors for the 
measured $l_j^m$. Consider the translational potential function
\begin{align}
\cU_t (\hat \msg,\bar a^m,\bar p) &= \frac12 \kappa y\T y= \frac12\kappa \|\bar{p}-\hat{R}\bar{a}^m-\hat{b}\|^2,
\label{U0t} 
\end{align}
where $\bar{p}$ is defined by \eqref{barp}, $y\equiv y(\hat\msg,\bar{a}^m,\bar p)=\bar{p}-
\hat{R}\bar{a}^m-\hat{b}$ and $\kappa$ is a positive scalar. Therefore, the total potential function is 
defined as the sum of the generalization of \eqref{U0r} defined in~\cite{Automatica,ast_acc14}  
for attitude determination on $\SO$, and the translational energy \eqref{U0t} as
\begin{align}
\cU(\hat{\msg},L^m,D,\bar a^m,\bar p)&= \Phi \big( \cU^0_r (\hat{\msg},L^m,D) \big)+\cU_t(\hat{\msg},\bar a^m,\bar p)\nn\\
&=\Phi \big(\frac12\lan D -\hat R L^m , (D -\hat R L^m)W\ran\big)\nn\\
&~~~~~~~+\frac12\kappa \|\bar{p}-\hat{R}\bar{a}^m-\hat{b}\|^2,\label{costU}
\end{align}
where $W$ is positive definite (not necessarily diagonal) which can be selected according to Lemma 3.2 in \cite{Automatica}, and 
$\Phi: [0,\infty)\mapsto[0,\infty)$ is a $\mC^2$ function that satisfies $\Phi(0)=0$ and 
$\Phi'(\mpz x)>0$ for all $\mpz x\in[0,\infty)$. Furthermore, $\Phi'(\cdot)\leq\alpha(\cdot)$ where 
$\alpha(\cdot)$ is a Class-$\mathcal{K}$ function~\cite{khal} and $\Phi'(\cdot)$ denotes the 
derivative of $\Phi(\cdot)$ with respect to its argument. Because of these properties 
of the function $\Phi$, the critical points and their indices coincide for $\cU^0_r$ and 
$\Phi(\cU^0_r)$~\cite{Automatica}. Define the kinetic energy-like function: 
\be
\cT \Big(\varphi(\hat\msg,\xi^m,\hat\xi)\Big)= \frac12 \varphi(\hat\msg,\xi^m,\hat\xi)\T \bJ\varphi(\hat\msg,\xi^m,\hat\xi), 
\label{costT} \ee
where $\bJ\in\bR^{6\times 6}>0$ is an artificial inertia-like kernel matrix. Note that in contrast to 
rigid body inertia matrix, $\bJ$ is not subject to intrinsic physical constraints like the triangle 
inequality, which dictates that the sum of any two eigenvalues of the inertia matrix has to be 
larger than the third. Instead, $\bJ$ is a gain matrix that can be used to tune the estimator. For 
notational convenience, $\varphi(\hat\msg,\xi^m,\hat\xi)$ is denoted as $\varphi$ from 
now on; this quantity is the velocities estimation error in the absence of measurement 
noise. Now define the Lagrangian 
\be
\cL (\hat{\msg},L^m,D,\bar a^m,\bar p,\varphi)= \cT(\varphi) -\cU(\hat{\msg},L^m,D,\bar a^m,\bar p),
\label{contLag}\ee
and the corresponding action functional over an arbitrary time interval $[t_0,T]$ for $T>0$,
\be \cS \big(\cL (\hat{\msg},L^m,D,\bar a^m,\bar p,\varphi)\big)= \int_{t_0}^T \cL 
(\hat{\msg},L^m,D,\bar a^m,\bar p,\varphi) \di t, \,\label{action} \ee
such that $\dot{\hat \msg}= \hat \msg(\hat\xi)^\vee$. A Rayleigh dissipation term linear in the velocities of the form $\bD\varphi$ 
where $\bD\in\bR^{6\times 6}>0$ is used in addition to the Lagrangian \eqref{contLag}, and 
the Lagrange-d'Alembert principle from variational mechanics is applied to obtain the estimator on 
$\Ta\SE$. This yields
\begin{align}
\delta_{\msh,\varphi} \cS \big(\cL (\msh,D,\bar{p},\varphi)\big)=\int_{t_0}^T \eta\T
\bD\varphi \di t,\label{LagdAlem}
\end{align}
which in turn results in the following continuous-time filter.

\subsection{Variational Estimator for Pose and Velocities}
The nonlinear variational estimator obtained by 
applying the Lagrange-d'Alembert principle to the Lagrangian \eqref{contLag} with a dissipation 
term linear in the velocities estimation error, is given by the following statement.

\begin{theorem}\label{filterTHM}
The nonlinear variational  estimator for pose and velocities is given by
\begin{align}
\begin{cases}
\bJ\dot{\varphi}&=\adast{\varphi}\bJ\varphi-Z(\hat{\msg},L^m,D,\bar{a}^m,\bar{p})-\bD\varphi,\vspace{.05in}\\
\hat{\xi}&=\xi^m-\Ad{\hat\msg^{-1}}\varphi,
\vspace{.05in}\\
\dot{\hat{\msg}}&=\hat{\msg}(\hat{\xi})^\vee,\label{ContFil}
\end{cases}
\end{align}
where $\adast{\zeta}=(\ad{\zeta})\T$ with $\ad{\zeta}$ defined by
\begin{align}
\ad{\mpz{\zeta}}=\bbm \mpz w^\times~~ & 0\\ \mpz v^\times\;\; & \mpz{w}^\times\ebm\mbox{ for } \zeta=\bbm \mpz w\\ \mpz v\ebm\in\bR^6,\label{ad_def}
\end{align}
and $Z(\hat{\msg},L^m,D,\bar{a}^m,\bar{p})$ is defined by
\begin{align}
\begin{split}
Z(\hat{\msg},L^m,D,&\bar{a}^m,\bar{p})=\label{Z}\\
&\bbm \Phi'\Big(\cU^0_r (\hat{\msg},L^m,D)\Big)S_\Gamma(\hat{R})+\kappa\bar{p}^\times y\\ \kappa y\ebm,\end{split}
\end{align}
where $\cU^0_r (\hat{\msg},L^m,D)$ is defined as \eqref{U0r}, $y\equiv y(\hat{\msg},\bar{a}^m,\bar{p})=\bar{p}-\hat{R}\bar{a}^m-\hat{b}$ and
\begin{align} S_\Gamma(\hat{R})=\mrm{vex}\big(DW(L^m)\T\hat{R}\T-\hat{R}L^mWD\T\big), \label{SLdef} \end{align}
where $\mrm{vex}(\cdot): \so\to\bR^3$ is the inverse of the $(\cdot)^\times$ map.
\end{theorem}
The proof is presented in \cite{Automatica2,Gaurav_ASR}. In the proposed approach, the time 
evolution of $(\hat \msg,\hat\xi)$ has the form of the dynamics of a rigid body with Rayleigh 
dissipation. This results in an estimator for the motion states $(\msg,\xi)$ that 
dissipates the ``energy" content in the estimation errors $(\mathsf{h},\varphi)= (\msg \hat\msg^{-1}, 
\Ad{\hat\msg}(\xi- \hat\xi))$ to provide guaranteed asymptotic stability in the case of perfect 
measurements~\cite{Automatica}. The variational pose estimator can also be interpreted as a low-pass stable filter (cf. \cite{Tayebi2011}). 
Indeed, one can connect the low-pass filter interpretation to the simple example of the natural dynamics 
of a mass-spring-damper system. This is a consequence of the fact that the mass-spring-damper system 
is a mechanical system with passive dissipation, evolving on a configuration space that is the vector 
space of real numbers, $\bR$. In fact, the equation of motion of this system can be obtained by 
application of the Lagrange-d'Alembert principle on the configuration space $\bR$. 
If this analogy or interpretation is extended to a system evolving on a Lie group as a 
configuration space, then the generalization of the mass-spring-damper system is a ``forced 
Euler-Poincar\'{e} system'' with passive dissipation, as is obtained here.

\section{DISCRETIZATION FOR COMPUTER IMPLEMENTATION}\label{Sec5}
For onboard computer implementation, the variational estimation scheme outlined above has to 
be discretized. Since the estimation scheme proposed here is obtained from 
a variational principle of mechanics, it can be discretized by applying the discrete 
Lagrange-d'Alembert principle~\cite{marswest}. Consider an interval of time $[t_0, T]\in\bR^+$ 
separated into $N$ equal-length subintervals $[t_i,t_{i+1}]$ for $i=0,1,\ldots,N$, with $t_N=T$ 
and $t_{i+1}-t_i=\Delta t$ is the time step size. Let $(\hat \msg_i,\hat\xi_i)\in\SE\times\bR^6$ 
denote the discrete state estimate at time $t_i$, such that $(\hat \msg_i,\hat\xi_i)\approx 
(\hat \msg(t_i),\hat\xi(t_i))$ where $(\hat \msg(t),\hat\xi(t))$ is the exact solution of the 
continuous-time estimator at time $t\in [t_0, T]$. Let the values of the discrete-time measurements  
$\xi^m$, $\bar a^m$ and $L^m$ at time $t_i$ be denoted as $\xi^m_i$, $\bar a^m_i$ and $L^m_i$, 
respectively. Further, denote the corresponding values for the latter two quantities in inertial frame 
at time $t_i$ by $\bar p_i$ and $D_i$, respectively. The discrete-time filter is then presented in the form of a Lie group variational integrator (LGVI) in the following statement.

\begin{theorem} \label{discfilter}
A first-order discretization of the estimator proposed in Theorem \ref{filterTHM} is given by
\begin{align}
(J\omega_i)^\times&=\frac{1}{\Delta t}(F_i\mathcal{J}-\mathcal{J}F_i\T),\label{LGVI_F}\\
(M+\Delta t\bD_t)\upsilon_{i+1}&=F_i\T M\upsilon_i\label{LGVI_upsilon}\\
&~~~~~~~~+\Delta t \kappa (\hat{b}_{i+1}+\hat{R}_{i+1}\bar{a}^m_{i+1}-\bar{p}_{i+1}),\nn\\
(J+\Delta t\bD_r)\omega_{i+1}&=F_i\T J\omega_i+\Delta t M\upsilon_{i+1}\times\upsilon_{i+1}\nn\\
+\Delta t&\kappa \bar{p}_{i+1}^\times (\hat{b}_{i+1}+\hat{R}_{i+1}\bar{a}^m_{i+1})\label{LGVI_omega}\\
-\Delta t&\Phi' \big( \cU^0_r (\hat{\msg}_{i+1},L_{i+1}^m,D_{i+1}) \big)S_{\Gamma_{i+1}}(\hat{R}_{i+1}),\nn\\
\hat\xi_i&=\xi^m_i-\Ad{\hat\msg_i^{-1}}\varphi_i,\label{LGVI_xihat}\\
\hat\msg_{i+1}&=\hat\msg_i\exp(\Delta t\hat\xi_i^\vee),\label{LGVI_ghat}
\end{align}
where $F_i\in\SO$, $\big(\hat\msg(t_0),\hat\xi(t_0)\big)=(\hat\msg_0,\hat\xi_0)$, $\mathcal{J}$ is defined in terms of positive
matrix $J$ by $\mathcal{J}=\frac12\tr[J]I-J$, $M$ is a positive definite matrix,
$\varphi_i=[\omega_i\T\;\upsilon_i\T]\T$, and $S_{\Gamma_i}(\hat R_i)$ is the value of  
$S_\Gamma(\hat R)$ at time $t_i$, with $S_\Gamma(\hat R)$ defined by \eqref{SLdef}.
\end{theorem}

\section{NUMERICAL SIMULATIONS}\label{Sec6}
\begin{figure}
\begin{center}
\includegraphics[height=3.8in]{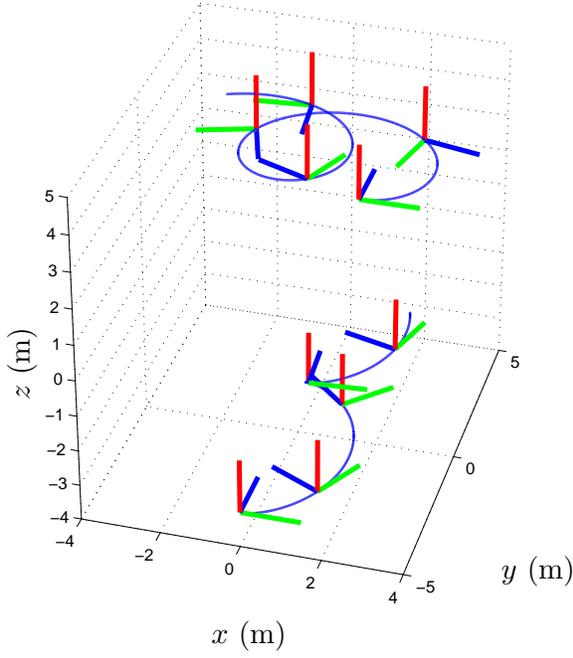}
\caption{Position and attitude trajectory of the simulated vehicles.}  
\label{Fig1}                                 
\end{center}                                 
\end{figure}

\begin{figure}
\begin{center}
\includegraphics[height=2.4in]{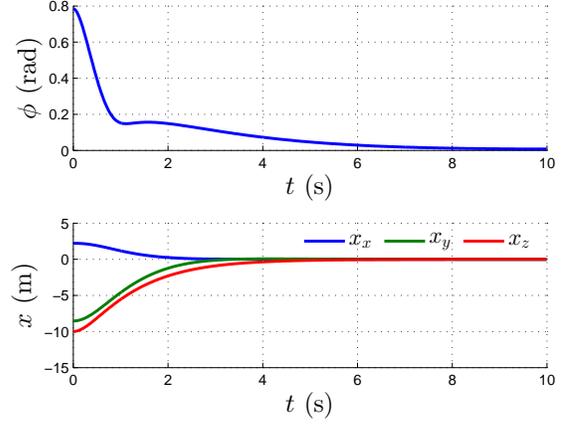}
\caption{Principal angle of the relative attitude and position estimation error.}  
\label{Fig2}                                 
\end{center}                                 
\end{figure}

\begin{figure}
\begin{center}
\includegraphics[height=2.4in]{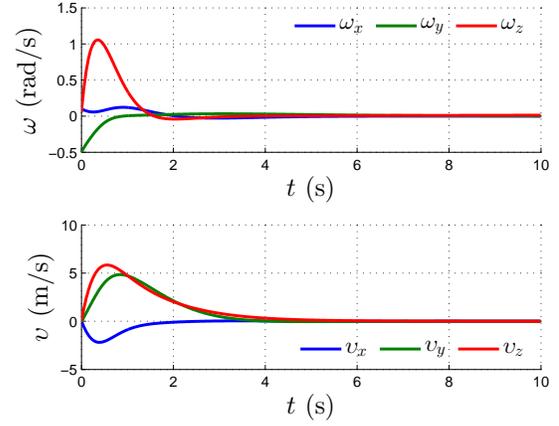}
\caption{Relative angular and translational velocity estimation error.}  
\label{Fig3}                                 
\end{center}                                 
\end{figure}

This section presents numerical simulation results of the discrete time estimator described in 
Section \ref{Sec5}, which is a Lie group variational integrator. Consider two vehicles performing 
spatial maneuvers, as shown in Fig.~\ref{Fig1}. These trajectories are generated using the 
equations of motion for these two vehicles and in turn generate the ``true" relative states of 
one vehicle with respect to another. The UAV at higher altitude has a camera that has the 
lower UAV in its FOV at all instants. 
The initial relative attitude and relative position of the lower vehicle with respect to the  
higher vehicle, are:
\begin{align}
R_0=I\; \mbox{ and } b_0=[1.5\;\;\;\; 5\;\;\;\; 6]\T\mbox{ m}.
\end{align}
The initial relative angular and relative translational velocity of these two vehicles are:
\begin{align}
\begin{split}
\omega_0=0\mbox{ rad/s}, \mbox{ and } \nu_0=[0.08\;\;-0.003\;\;-0.0007]\T\mbox{ m/s}.
\end{split}
\end{align}
There are three feature points on the lower vehicle's body, and their positions expressed in 
the lower vehicle's body frame are
\begin{align}
P=\bbm1\;\; &0\;\; &0\\0 &1 &-1\\0 &0 &0\ebm.
\end{align}
Relative position vectors of these points are measured by the camera on the upper vehicle. 
Velocities of these points are calculated using the linear filter introduced in \cite{ACC2015}. 
The relative velocities can be computed using these measurements by \eqref{ximMore2}. All 
the camera readings contain random zero mean signals whose probability distributions are 
normalized bump functions with the width equal to 1 mm in each coordinate. The ``inertia-like" 
gain matrices for the estimator are selected to be: 
\begin{align}
\begin{split}
J&=\diag\big([0.9\;\;0.6\;\;0.3]\T\big),\\
M&=\diag\big([0.0608\;\;0.0486\;\;0.0365]\T\big)
\end{split}
\end{align}
The ``dissipation" gain matrices for the estimator are set to: 
\begin{align}
\begin{split}
\bD_r&=\diag\big([2.7\;2.2\;1.5]\T\big),\bD_t=\diag\big([0.1\;\;0.12\;\;0.14]\T\big).
\end{split}
\end{align}
$\Phi(\cdot)$ could be any $C^2$ function with the properties described in Section \ref{Sec3}, but is 
selected to be $\Phi(x)=x$ here. The initial estimated states have the following values:
\begin{align}
\hat R_0&=\expm_{\SO}\big((\frac{\pi}{4}\times[0\ 0\ 1]\T)^\times
\big),\;\hat b_0=[-3\;\;\;2\;\;\;4]\T\mbox{ m}\nn\\
\hat\omega_0&=[0.1\;\;-0.5\;\;\;0.05]\T\mbox{ rad/s},\\
\mbox{ an}&\mbox{d }\hat\nu_0=[0.05\;\;-0.09\;\;\;0.01]\T\mbox{ m/s}.\nn
\end{align}
The discrete-time estimator \eqref{LGVI_F}-\eqref{LGVI_ghat} is simulated over a time interval of 
$T=10$ s with time stepsize $h=0.01$ s. At each instant, \eqref{LGVI_F} is solved using 
Newton-Raphson iterations to find $F_i$. Then, the rest of the equations (all explicit) are solved 
consecutively to generate the estimated states. The principal angle of the relative attitude 
estimation error and components of the relative position estimate error are plotted in 
Fig.~\ref{Fig2}. Components of the relative angular and translational velocities are depicted 
in Fig.~\ref{Fig3}.

As can be noticed from the figures, all the estimated relative states converge to a bounded 
neighborhood of the corresponding true relative states, where the size of this neighborhood 
depends on the level of measurement noise and estimator gains. This confirms the stability 
and convergence properties of the estimator.

\section{CONCLUSION}\label{Sec7}
This article proposes an estimator for relative pose and relative velocities of one vehicle with 
respect to another vehicle that uses only optical measurements from onboard optical sensor(s). The 
sensors are assumed to provide measurements in continuous-time or at a high frequency, with 
bounded measurement noise due to limited fields of view. A Lagrangian in terms of measurement  
residuals and which can be expressed in terms of state estimation errors when perfect 
measurements are available, is proposed. Applying the Lagrange-d'Alembert principle to this 
Lagrangian with a dissipation term linear in relative velocity estimation errors, an estimator is 
designed on the Lie group of relative motions between two rigid vehicles. In the case of perfect 
measurements, this estimator is shown to be almost globally asymptotically stable with a domain 
of convergence that is open and dense in the state space. The continuous estimator is discretized 
by applying the discrete Lagrange-d'Alembert principle to the discretized Lagrangian and 
dissipation terms for rotational and translational motions. In the presence of measurement noise, 
numerical simulations with this discrete estimator show that state estimates converge to a 
bounded neighborhood of the true relative motion states. Future work will be directed towards 
creating higher-order discretizations of the continuous-time filter given here.  
\pagestyle{empty}

\addtolength{\textheight}{-2cm}   



%
%
%



\begin{thebibliography}{10}

\bibitem{aguhes06}
Aguiar,~A., \& Hespanha,~J. (2006).
\newblock Minimum-energy state estimation for systems with perspective outputs.
\newblock {\em IEEE Transactions on Automatic Control}, 51(2), 226--241.

\bibitem{amelin2014algorithm}
Amelin,~K.~S., \& Miller,~A.~B. (2014).
\newblock An algorithm for refinement of the position of a light UAV on the
  basis of Kalman filtering of bearing measurements.
\newblock {\em Journal of Communications Technology and Electronics},
  59(6), 622--631.

\bibitem{Bayadi2014almost}
Bayadi,~R., \& Banavar,~R.~N. (2014).
\newblock Almost global attitude stabilization of a rigid body for both internal and external actuation schemes.
\newblock {\em European Journal of Control}, 20(1), 45--54.

%

\bibitem{bonmaro09}
Bonnabel,~S., Martin,~P., \& Rouchon,~P. (2009).
\newblock Nonlinear symmetry-preserving observers on {L}ie groups.
\newblock {\em IEEE Transactions on Automatic Control}, 54(7), 1709--1713.

\bibitem{chaturvedi2011rigid}
Chaturvedi,~N.~A., Sanyal,~A.~K., \& McClamroch,~N.~H. (2011).
\newblock Rigid-body attitude control.
\newblock {\em IEEE Control Systems Magazine}, 31(3), 30--51.

\bibitem{Good2013ECC}
Goodarzi,~F.~A., Lee,~D., \& Lee,~T. (2013).
\newblock Geometric nonlinear {PID} control of a quadrotor {UAV} on {SE}(3).
\newblock In {\em Proceedings of the European Control Conference} (pp. 3845--3850). Zurich, Switzerland.

\bibitem{Good2014ASME}
Goodarzi,~F.~A., Lee,~D., \& Lee,~T. (2014).
\newblock Geometric Adaptive Tracking Control of a Quadrotor UAV on SE(3) for Agile Maneuvers.
\newblock {\em ASME Journal of Dynamic Systems, Measurement and Control}, 137(9), 091007.

\bibitem{Good2014ACC}
Goodarzi,~F.~A., Lee,~D., \& Lee,~T. (2014).
\newblock Geometric stabilization of a quadrotor UAV with a payload connected by flexible cable.
\newblock In {\em Proceedings of the American Control Conference} (pp. 4925--4930). Portland, OR, USA.



\bibitem{ICRA2015}
Izadi,~M., Samiei,~E., Sanyal,~A.~K., \& Kumar,~V. (2015).
\newblock Comparison of an attitude estimator based on the
  {L}agrange-d'{A}lembert principle with some state-of-the-art filters.
\newblock In {\em Proceedings of the IEEE International Conference on Robotics and
  Automation} (pp. 2848--2853). Seattle, WA, USA.

\bibitem{ACC2015}
Izadi,~M., Sanyal,~A.~K., Samiei,~E., \& Viswanathan,~S.~P. (2015).
\newblock Discrete-time rigid body attitude state estimation based on the
  discrete {L}agrange-d'{A}lembert principle.
\newblock In {\em Proceedings of the American Control Conference} (pp. 3392--3397). Chicago, IL, USA.

\bibitem{Automatica}
Izadi,~M., \& Sanyal,~A.~K. (2014).
\newblock Rigid body attitude estimation based on the {L}agrange-d'{A}lembert
  principle.
\newblock {\em Automatica}, 50(10), 2570--2577.

\bibitem{Automatica2}
Izadi,~M., \& Sanyal,~A.~K. (2015).
\newblock Rigid body pose estimation based on the {L}agrange-d'{A}lembert principle.
\newblock To appear in {\em Automatica}.

\bibitem{khal}
Khalil,~H.~K. (2001).
\newblock {\em Nonlinear Systems} (3$^{rd}$ edition).
\newblock Prentice Hall, Upper Saddle River, NJ.

\bibitem{Khosravian2015Recursive}
Khosravian,~A., Trumpf,~J., Mahony,~R., \& Hamel,~T. (2015).
\newblock Recursive Attitude Estimation in the Presence of Multi-rate and Multi-delay Vector Measurements.
\newblock In {\em Proceedings of the American Control Conference} (pp. 3199--3205). Chicago,
  IL, USA.

\bibitem{Khosravian2015observers}
Khosravian,~A., Trumpf,~J., Mahony,~R., \& Lageman,~C. (2015).
\newblock Observers for invariant systems on Lie groups with biased input
  measurements and homogeneous outputs.
\newblock {\em Automatica}, 55, 19--26.



\bibitem{leishman2014relative}
Leishman,~R.~C., McLain,~T.~W., \& Beard,~R.~W. (2014).
\newblock Relative navigation approach for vision-based aerial {GPS}-denied
  navigation.
\newblock {\em Journal of Intelligent \& Robotic Systems}, 74(1-2), 97--111.

\bibitem{mahapf08}
Mahony,~R., Hamel,~T., \& Pflimlin,~J.~M. (2008).
\newblock Nonlinear complementary filters on the special orthogonal group.
\newblock {\em IEEE Transactions on Automatic Control}, 53(5), 1203--1218.

\bibitem{mabeda04}
Maithripala,~D.~H., Berg,~J.~M., \& Dayawansa,~W.~P. (2004).
\newblock An intrinsic observer for a class of simple mechanical systems on a
  {L}ie group.
\newblock In {\em Proceedings of the American Control Conference} (pp. 1546--1551). Boston,
  MA, USA.

\bibitem{markSO3}
Markley,~F.~L. (2006).
\newblock Attitude filtering on {SO}(3).
\newblock {\em The Journal of the Astronautical Sciences}, 54(4), 391--413.


\bibitem{marswest}
Marsden,~J.~E., \& West,~M. (2001).
\newblock Discrete mechanics and variational integrators.
\newblock {\em Acta Numerica}, 10, 357--514.

\bibitem{miller2014tracking}
Miller,~A., \& Miller,~B. (2014).
\newblock Tracking of the UAV trajectory on the basis of bearing-only
  observations.
\newblock In {\em Proceedings of the 53rd Annual Conference on Decision and Control} (pp. 4178--4184). Los Angeles, CA, USA.


\bibitem{Gaurav_ASR}
Misra,~G., Izadi,~M., Sanyal,~A.~K., \& Scheeres,~D.~J. (2015).
\newblock Coupled orbit-attitude dynamics and relative state estimation of
  spacecraft near small Solar System bodies.
\newblock {\em Advances in Space Research}.

\bibitem{Mortensen}
Mortensen,~R.~E. (1968).
\newblock Maximum-likelihood recursive nonlinear filtering.
\newblock {\em Journal of Optimization Theory and Applications}, 2(6), 386--394.

\bibitem{rehbinder2003pose}
Rehbinder,~H., \& Ghosh,~B.~K. (2003).
\newblock Pose estimation using line-based dynamic vision and inertial sensors.
\newblock {\em IEEE Transactions on Automatic Control}, 48(2), 186--199.

\bibitem{sanyal2009inertia}
Sanyal,~A.~K., Fosbury,~A., Chaturvedi,~N.~A., \& Bernstein,~D.~S. (2009).
\newblock Inertia-free spacecraft attitude tracking with disturbance rejection
  and almost global stabilization.
\newblock {\em Journal of Guidance, Control, and Dynamics}, 32(4), 1167--1178.

\bibitem{ast_acc14}
Sanyal,~A.~K., Izadi,~M., \& Butcher,~E.~A. (2014).
\newblock Determination of relative motion of a space object from simultaneous
  measurements of range and range rate.
\newblock In {\em Proceedings of the American Control Conference} (pp. 1607--1612). Portland, OR, USA.

\bibitem{sany}
Sanyal,~A.~K., Lee,~T., Leok,~M., \& McClamroch,~N.~H. (2008).
\newblock Global optimal attitude estimation using uncertainty ellipsoids.
\newblock {\em Systems \& Control Letters}, 57(3), 236--245.



\bibitem{shen2013vision}
Shen,~S., Mulgaonkar,~Y., Michael,~N., \& Kumar,~V. (2013).
\newblock Vision-based state estimation and trajectory control towards
  aggressive flight with a quadrotor.
\newblock In {\em Proceedings of the Robotics Science and Systems}.

\bibitem{shen2013rotor}
Shen,~S., Mulgaonkar,~Y., Michael,~N., \& Kumar,~V. (2013).
\newblock Vision-based state estimation for autonomous rotorcraft {MAV}s in
  complex environments.
\newblock In {\em Proceedings of the IEEE International Conference on Robotics and
  Automation} (pp. 1758--1764). Karlsruhe, Germany.

\bibitem{Tayebi2011}
Tayebi,~A., Roberts,~A., \& Benallegue,~A. (2011).
\newblock Inertial measurements based dynamic attitude estimation and velocity-free attitude stabilization.
\newblock In {\em Proceedings of the American Control Conference} (pp. 1027--1032). San Francisco, CA, USA.

\bibitem{Vas1}
Vasconcelos,~J.~F., Cunha,~R., Silvestre,~C., \& Oliveira,~P. (2010).
\newblock A nonlinear position and attitude observer on {SE}(3) using landmark
  measurements.
\newblock {\em Systems \& Control Letters}, 59, 155--166.

\bibitem{silvest08}
Vasconcelos,~J.~F., Silvestre,~C., \& Oliveira,~P. (2008).
\newblock A nonlinear {GPS/IMU} based observer for rigid body attitude and
  position estimation.
\newblock In {\em Proceedings of the IEEE Conference on Decision and Control} (pp. 1255--1260). Cancun, Mexico.

\bibitem{jo:wahba}
Wahba,~G. (1965).
\newblock A least squares estimate of satellite attitude, Problem 65-1.
\newblock {\em SIAM Review}, 7(5), 409.

\bibitem{ZamPhD}
Zamani,~M. (2013).
\newblock {\em Deterministic Attitude and Pose Filtering, an Embedded {L}ie
  Groups Approach}.
\newblock Ph.D. Thesis. Australian National University, Canberra, Australia.

\bibitem{Zamani2013minimum}
Zamani,~M., Trumpf,~J., \& Mahony,~R. (2013).
\newblock Minimum-energy filtering for attitude estimation.
\newblock {\em IEEE Transactions on Automatic Control}, 58(11), 2917--2921.

\end{thebibliography}
\end{document}